\documentclass[11pt,A4]{style}

\usepackage{graphicx}
\usepackage{amsmath}
\usepackage{amsfonts}
\usepackage{amssymb}

\usepackage{hyperref}
\usepackage{amsthm}

\newtheorem{theorem}{Theorem}[section]
\newtheorem{Proposition}{Proposition}[section]

\newtheorem{Corollary}[theorem]{Corollary}

\newtheorem{Problem}[theorem]{Problem}

\theoremstyle{definition} %

\theoremstyle{remark} %

\usepackage{xcolor}

\newcommand{\bs}[1]{{#1}} 

\newcommand{\hOmega}{\widehat\Omega}

\newcommand{\hGamma}{\widehat\Gamma}
\newcommand{\hI}{\widehat I}

\newcommand{\hX}{\bs{\widehat{X}}}

\newcommand{\hV}{\bs{\widehat{V}}}

\newcommand{\hL}{{\widehat{L}}}

\newcommand{\hF}{\bs{\widehat{F}}}
\newcommand{\hJ}{\widehat{J}}
\newcommand{\hE}{\bs{\widehat{E}}}

\newcommand{\hdiv}{\widehat{\operatorname{div}}}

\newcommand{\hnabla}{\bs{\widehat\nabla}}
\newcommand{\hSigma}{\bs{\widehat\Sigma}}

\newcommand{\hsigma}{\bs{\widehat\sigma}}

\newcommand{\hrho}{\bs{\widehat\rho}}

\newcommand{\hs}{\bs{\widehat s}}

\newcommand{\hd}{\bs{\widehat d}}

\newcommand{\hv}{\bs{\widehat v}}
\newcommand{\hu}{\bs{\widehat u}}
\newcommand{\hp}{\bs{\widehat p}}

\newcommand{\hf}{\bs{\widehat f}}
\newcommand{\hg}{\bs{\widehat g}}

\newcommand{\hA}{\bs{\widehat A}}

\newcommand{\hU}{\bs{\widehat U}}

\newcommand{\hz}{\bs{\widehat z}}
\newcommand{\hZ}{\bs{\widehat Z}}

\newcommand{\hPsi}{\bs{\widehat\Psi}}
\newcommand{\hPhi}{\bs{\widehat\Phi}}
\newcommand{\hpsi}{\bs{\widehat\psi}}

\newcommand{\hn}{\bs{\widehat n}}

\newcommand{\calL}{{\cal L}}

\title{Adjoint-based methods for optimization and goal-oriented error control applied to fluid-structure interaction:
implementation of a partition-of-unity dual-weighted residual estimator for stationary forward FSI problems in deal.II}

\author{Thomas Wick}

\address{Leibniz University Hannover\\
Institute of Applied Mathematics\\
Hannover, Germany\\
e-mail: thomas.wick@ifam.uni-hannover.de
\and
Cluster of Excellence PhoenixD\\ (Photonics, Optics, and
	Engineering - Innovation Across Disciplines)\\ 
Leibniz University Hannover, Germany
}

\keywords{goal-oriented error control, dual-weighted residuals, adjoint, 
mesh adaptivity, fluid-structure interaction, deal.II}

\abstract{In this work, we implement 
goal-oriented error control and spatial mesh adaptivity 
for stationary fluid-structure interaction. The a posteriori 
error estimator is realized using the dual-weighted residual method
in which the adjoint equation arises. The fluid-structure interaction
problem is formulated within a variational-monolithic framework
using arbitrary Lagrangian-Eulerian coordinates. The overall 
problem is nonlinear and solved with Newton's method. 
We specifically consider the FSI-1 benchmark problem 
in which quantities of interest include 
the elastic beam displacements, drag, and lift.
The implementation is provided
open-source published on github \texttt{https://github.com/tommeswick/goal-oriented-fsi}.
Possible extensions are discussed in the source code 
and in the conclusions of this paper.
}

\begin{document}

\section{Introduction}
Fluid-structure interaction (FSI) is well-known \cite{BuSc06,FoQuaVe09,GaRa10,BuSc10,BaTaTe13,BoGaNe14,Ri17_fsi,FrHoRiWiYa17,Wi20_book} and a prime example 
of a multiphysics problem. It combines several challenges such as 
different types of partial-differential equations (PDE), interface-coupling, 
nonlinearities in the equations and due to coupling, Lagrangian and Eulerian coordinates.
These result into typical numerical challenges such as robust spatial discretization (in particular 
for the moving interface), robust time-stepping schemes, efficient and robust linear 
and nonlinear solution algorithms. Computational work include different coupling 
concepts \cite{HuLiZi81,Glowinski2001363,TeSaKeSt06,Pe02,Du06},
space-time multiscale \cite{TaTe11}, 
reduced order modeling \cite{FoGeNoQua01,Lassila2013,TeCoBai20,hagmeyer2021fluidbeam}, 
optimal control, parameter estimation, uncertainty quantification 
\cite{PeVeVe11,BeMoiGer12,RiWi13_fsi_opt,Kratzke18,FaiRi20,WiWo20}, 
and efficient solver developments 
\cite{He04,Barker2010642,gee2011truly,RazzaqDamanikHronOuazziTurek:2012,CrDeFouQua11,Ri15_fsi_solver,DeparisFortiGrandperrinQuarteroni:2016,JoLaWi19_fsi,Wich2021}.

In this work, the main objective is the 
application and open-source implementation of goal-oriented 
a posteriori error control using the dual-weighted residual (DWR)
method \cite{BeRa01,BaRa03}. For applications 
in fluid-structure interaction, we refer to 
\cite{GraeBa06,FickBrummelenZee2010,VANDERZEE20112738,Wi11_dwr_Nov_2011,Ri12_dwr,Ri17_fsi,Fai17,FaiWi18}.
A recent overview of our own work using the adjoint FSI equation in goal-oriented 
error estimation and optimization was done in \cite{Wi21_WCCM}.
In \cite{RiWi15_dwr} a variational localization using a partition-of-unity (PU) was proposed,
facilitating the application to coupled problems such as fluid-structure interaction.
In view of increasing initiatives of open-source developments, another purpose 
of this work is to provide a documented open-source code. To this end, 
a stationary fluid-structure interaction problem is considered in order to 
explain the main steps of a PU-DWR estimator. The problem is formulated 
within a monolithic framework using arbitrary Lagrangian-Eulerian (ALE) coordinates.
For some well-posedness results of such stationary FSI problems, 
we refer to \cite{Gr02,WiWo19}.
Together with the goal functional under consideration, 
the FSI formulation serves as PDE constraint and the Lagrange formalism 
can be applied. Specifically, the monolithic formulation yields a consistent 
adjoint equation. 

For the monolithic, stationary, FSI formulation we follow \cite{RiWi10,Wi11_phd}
and for the PU-DWR error estimator, we follow \cite{RiWi15_dwr}. 
The basis of our programming code is \cite{Wi13_fsi_with_deal} 
(see also updates on github\footnote{\url{https://github.com/tommeswick/fsi}}) and we take some ideas from
deal.II \cite{dealII91,deal2020} 
step-14\footnote{\url{https://www.dealii.org/current/doxygen/deal.II/step_14.html}}.
Our resulting code can be found on github\footnote{\url{https://github.com/tommeswick/goal-oriented-fsi}}.

\section{Variational-monolithic ALE fluid-structure interaction}

\subsection{Modeling}
For the function spaces in the (fixed) reference domains 
$\hOmega, \hOmega_f,\hOmega_s$, we define 
$
\hV:= H^1 (\hOmega)^d.
$
In the fluid and solid domains, we define further:
\begin{align*}
\hL_f   &:= L^2 (\hOmega_f), \quad
\hL_f^0 := L^2 (\hOmega_f)/\mathbb{R}, \quad
\hV_f^0 := \{ \hv_f\in H^1 (\hOmega_f)^d : \, 
\hv_f = 0 \text{ on } \hGamma_{\text{in}}\cup \hGamma_{D} \},\\
\hV_{f,\hu}^0 &:= \{ \hu_f\in H^1 (\hOmega_f)^d : \, 
\hu_f = \hu_s \text{ on } \hGamma_i, \quad \hu_f = 0 \text{ on } \hGamma_{\text{in}}
\cup \hGamma_{D}\cup\hGamma_{\text{out}} \}, \\
\hV_{f,\hu,\hGamma_i}^0 &:= \{ \hpsi_f\in H^1 (\hOmega_f)^d : 
\, \hpsi_f = 0 \text{ on } 
\hGamma_i \cup \hGamma_{\text{in}} \cup \hGamma_{D}\cup\hGamma_{\text{out}}\},\\
\hV_s^0 &:= \{ \hu_s\in H^1 (\hOmega_s)^d : \, 
\hu_s = 0 \text{ on } \hGamma_{D} \}.
\end{align*}
As stationary FSI problem in variational-monolithic ALE form, we have
\cite{Wi11_phd}[p. 29]:
\begin{Problem}
\label{fsi_equations_in_ale_stationary_problem}
Find $\{\hv_f,\hu_f,\hu_s,\hp_f\} 
\in \{ \hv_f^D + \hV_{f,\hv}^0\}\times 
\{ \hu_f^D + \hV_{f,\hu}^0\}\times \{ \hu_s^D + \hV_s^0 \}
\times\hL_f^0 $, such that
  \begin{eqnarray*}
    \begin{aligned}
      (\hrho_f \hJ  (\hF^{-1}\hv_f\cdot\hnabla) \hv_f),
      \hpsi^v)_{\hOmega_f} &\\
      + (\hJ\hsigma_f\hF^{-T},\hnabla\hpsi^v)_{\hOmega_f}      
      - \langle \hg_f, \hpsi^v \rangle_{\hGamma_N}
      - (\hrho_f \hJ\hf_f, \hpsi^v)_{\hOmega_f}
      &=0&&\forall\hpsi^v\in \hV_{f,\hv}^0, \\
      (\hF\hSigma ,\hnabla\hpsi^v)_{\hOmega_s}
      - (\hrho_s\hf_s, \hpsi^v)_{\hOmega_s}
      &=0&&\forall\hpsi^v\in \hV_s^0, \\    
      (\hsigma_{\text{mesh}} ,\hnabla\hpsi^u)_{\hOmega_f} +
      (\hv_s,\hpsi^u)_{\hOmega_s} 
      &=0&&\forall\hpsi^u\in \hV_{f,\hu,\hGamma_i}^0,\\
      (\hdiv\,(\hJ\hF^{-1}\hv_f),\hpsi^p)_{\hOmega_f}
       &=0&&\forall\hpsi^p\in \hL_f^0,
    \end{aligned}
  \end{eqnarray*}  
with $\hF = \hI + \hnabla\hu, \hJ = det(\hF), 
\hsigma_f = -\hp_f \hI + \hrho_f\nu_f (\hnabla\hv_f \hF^{-1} + \hF^{-T}\hnabla\hv_f), 
\hSigma = 2\mu_s \hE + \lambda_s tr(\hE)\hI, \hE = 0.5 (\hF^T\hF - \hI), \hsigma_{\text{mesh}} = \alpha_u \hnabla\hu_f$,
volume forces $\hf_f$ and $\hf_s$ (both zero in this work), flow correction term $\hg_f$ (do-nothing \cite{HeRaTu96}), densities $\hrho_s,\hrho_f$, kinematic viscosity $\nu_f$, and the Lam\'e parameters $\mu_s,\lambda_s$.
All explanations are provided in \cite{Wi11_phd}[Chapter 3].
\end{Problem}

\subsection{Discretization and numerical solution}
\label{sec_dis_primal}
For spatial discretization, a conforming Galerkin finite element scheme on 
quadrilateral mesh elements is employed \cite{Cia87}. Specifically, 
we use $Q_2^c$ elements for $\hv$ and $\hu:=\hu_f+\hu_s$, and $Q_1^c$ elements 
for $\hp$. For the flow problem $(\hv,\hp)$, 
this is the well-known inf-sup stable 
Taylor-Hood element; see e.g., \cite{GiRa1986}. 
Due to variational-monolithic coupling and globally-defined finite elements, the fluid pressure must be extended to the solid domain, which is achieved via $\alpha_u [(\hnabla\hp_s,\hnabla\hpsi^p) + (\hp_s,\hpsi^p)]$, and $\alpha_u$ (as before) small, positive. This 
is only for convenience, an alternative is to work with the 
\verb|FE_NOTHING|\footnote{\url{https://www.dealii.org/current/doxygen/deal.II/step_46.html}} element in deal.II.
The nonlinear 
problem is solved with Newton's method. Therein, for simplicity in this work,
we utilize a sparse direct solver \cite{DaDu97}. For algorithmic 
descriptions of our implementation, we refer to \cite{Wi11_phd}.

\newpage
\section{PU-DWR goal-oriented error control}
The Galerkin approximation reads:
Find $\hU_h=\{\hv_{f,h}, \hu_{f,h}, \hu_{s,h},\hp_{f,h}\}
\in \hX^0_{h,D}$, 
where $\hX_{h,D}^0 :=\{ \hv_{f,h}^D + \hV_{f,\hv,h}^0\} 
\times \{ \hu_{f,h}^D + \hV_{f,\hu,h}^0 \} \times \{ \hu_{s,h}^D + \hV_{s,h}^0 \}
\times \hL_{f,h}^0$, such that 
\begin{equation}\label{dis_abstract_stat_problem_discrete}
\hat A(\hU_h)(\hPsi_h ) = 
\hat F(\hPsi_h) \quad\forall \hPsi_h \in \hX_h,
\end{equation} 
where $\hX_h$ is the test space with homogeneous Dirichlet conditions.

\subsection{Goal functional}
The solution $\hU_h$ is used to calculate an approximation
$J(\hU_h)$ of the goal-functional $J(\hU):\hX\rightarrow\mathbb{R}$. 
This functional is assumed to be sufficiently differentiable.
The drag value as goal functional reads
\begin{equation*}
\label{aposteriori_line_integration_functional}
J(\hU) := \int_{\hat S} \hJ \hsigma_f \hF^{-T} \hn_f \, \hd \, \text{d}\hs,
\end{equation*}
where $\hn_f$ is the outward point normal vector of the cylinder boundary $\hat S$ \cite{HrTu06b} and the FSI interface $\hGamma_i$. Moreover, $\hd$ is a unit vector perpendicular to the mean flow direction. 
For the drag, we use $\hd = (1,0)$.

\subsection{Error representation}
We use the (formal) Euler-Lagrange method, 
to derive a computable representation of the approximation
error $J(\hU) - J(\hU_h)$. The task is: Find $\hU\in\hX^0_D$ such that
\begin{equation*}
\min \{J(\hU) - J(\hU_h)\} \quad\text{s.t. } \hA(\hU)(\hPsi) = 
\hat F(\hPsi) \quad\forall \hPsi \in \hX,
\end{equation*}
from which we obtain the optimality system 
\begin{align*}
\calL'_{\hZ} (\hU , \hZ)(\delta\hZ) &= \hF(\delta\hZ) 
- \hA(\hU)(\delta\hZ) = 0
\quad\forall\delta\hZ\in \hX, \quad\text{(Primal problem)}, \\
\calL'_{\hU} (\hU , \hZ)(\delta\hU) 
&= J'(\hU)(\delta\hU) - \hA'_{\hU}(\hU)(\delta\hU, \hZ) 
= 0 \quad\forall\delta\hU\in \hX, \quad\text{(Adjoint problem)}.
\end{align*}

Using the main theorem from \cite{BeRa01}, we obtain:
\begin{theorem}
We have the error identity:
\begin{align} \label{dwr_error_representation}
J(\hU) - J(\hU_h) = \frac{1}{2}\rho(\hU_h)(\hZ-\hPhi_h) +
\frac{1}{2}\rho^*(\hU_h, \hZ_h)(\hU-\hPsi_h) + {\cal R}^{(3)}_h,
\end{align}
for all $\{\hPsi_h , \hPhi_h\}\in \hX_h \times \hX_h$ and with the 
primal and adjoint residuals:
\begin{align*}
\rho(\hU_h)(\hZ-\hPhi_h) &:= -A(\hU_h)(\cdot) + \hat F(\cdot) , \\
\rho^*(\hU_h, \hZ_h)(\hU-\hPsi_h) &:= J'(\hU_h)(\cdot) 
- A'(\hU_h)(\cdot, \hZ_h) + \hat F(\cdot).
\end{align*}
The remainder term is ${\cal R}^{(3)}_h$ is of cubic order.
This error identity can be used to define the error estimator $\eta$, 
which can be further utilized to design adaptive schemes.
\end{theorem}
\begin{Corollary}[Primal error]
\label{coro_1}
The primal error identity reads:
\begin{align} \label{dwr_error_representation_primal}
J(\hU) - J(\hU_h) = \rho(\hU_h)(\hZ-\hPhi_h) + {\cal R}^{(2)}_h.
\end{align}
\end{Corollary}

\subsection{Adjoint equation, discretization, and numerical solution}
The adjoint equation reads: Find $\hZ = (\hz^v,\hz^u,\hz^p)\in\hX$ such that
\[
J'(\hU)(\hPhi) = \hA'_{\hU}(\hU)(\hPhi, \hZ) \quad\forall \hPhi\in\hX,
\]
and the explicit form can be found in \cite{Wi11_phd,Wi21_WCCM}.

For the discretization, we briefly mention that higher-order information 
for the adjoint solution must be employed due to Galerkin orthogonality;
in this work $\hX_h\subset \hX_h^{(2)} \subset \hX$.
For simplicity, this is realized with global-higher order finite 
elements and in order to ensure again inf-sup stability, we use 
$Q_4^c$ elements for $\hz^v$ and $\hz^u$, and $Q_2^c$ elements for $\hz^p$. 
It is clear that this is an expensive choice. For the numerical solution,
the same solvers as for the primal problem are taken (see Section \ref{sec_dis_primal}), 
namely a Newton-type method and sparse direct solver. Since the adjoint problem
is linear, Newton's method converges in one step. This is a trivial information, 
but for debugging reasons useful.

\subsection{Localization}
A PU localization \cite{RiWi15_dwr} for stationary FSI reads:
\begin{Proposition}
We have for the primal error part $\rho(\hU_h)(\cdot)$ 
the a posteriori error estimate
 \begin{equation}
\label{eq_est}
|J(\hU) - J(\hU_h)| \leq \eta := \bigl| \sum_{i=1}^M \eta_i \bigr| \leq \sum_{i=1}^M |\eta_i|
 \end{equation}
where $M$ is the dimension of the PU 
finite element space $\hV_{PU}$ (composed of $Q_1^c$ functions $\chi_i$)
and with the PU-DoF indicators
\begin{align*}
\eta_i &= - A(\hU_h)((\hZ_h^{(2)}-i_h \hZ_h^{(2)})\hPsi_i) 
+ \hat F((\hZ_h^{(2)}-i_h \hZ_h^{(2)})\hPsi_i)\\
&= -
(\hrho_f \hJ  (\hF^{-1}\hv_f\cdot\hnabla) \hv_f),\hpsi_i^v)_{\hOmega_f} 
- (\hJ\hsigma_f\hF^{-T},\hnabla\hpsi_i^v)_{\hOmega_f}      
+ \langle \hg_f, \hpsi_i^v \rangle_{\hGamma_N}\\
&\quad - (\hF\hSigma ,\hnabla\hpsi_i^v)_{\hOmega_s}
- (\hsigma_{\text{mesh}} ,\hnabla\hpsi_i^u)_{\hOmega_f}
- (\hdiv\,(\hJ\hF^{-1}\hv_f),\hpsi_i^p)_{\hOmega_f}\\
&\quad + (\hrho_f \hJ\hf_f,\hpsi_i^v)_{\hOmega_f} + (\hrho_s \hf_s,\hpsi_i^v)_{\hOmega_s}
\end{align*}
with the interpolation $i_h:\hX_h^{(2)}\to\hX_h$ and
the \textit{weighting functions} are defined as 
\begin{align*}
\hpsi_i^v &:= (\phi_{2h,v}^{(2)}-\phi_{h,v})\chi_i, \quad
\hpsi_i^u := (\phi_{2h,u}^{(2)}-\phi_{h,u})\chi_i, \quad
\hpsi_i^p := (\phi_{2h,p}^{(2)}-\phi_{h,p})\chi_i.
\end{align*}
\end{Proposition}

\subsection{Adaptive algorithm}
\begin{enumerate}
\item Compute the primal solution $\hU_h$ and the 
(higher-order) adjoint solution $\hZ_h^{(2)}$ 
on the present mesh $\mathcal{T}_h$.
\item Evaluate $|\eta| := |\sum_{i} \eta_{i}|$ in \eqref{eq_est}.
\item Check, if the stopping criterion is satisfied:
$|J(\hU) - J(\hU_h)| \leq |\eta| \leq TOL$, then accept $U_h$ 
within the tolerance $TOL$. 
Otherwise, proceed
to the following step. 
\item Mark all elements $K_i$ for refinement that touch DoFs $i$ 
with indicator $\eta_{i}$ with 
$\eta_{i} \geq \frac{\alpha\eta}{M_{el}}$ (where
$M_{el}$ denotes the total number of elements of the mesh $\mathcal{T}_h$
and $\alpha \approx 1$). \\
Alternatively, pure DoF-based refinement in $i$ can be carried out.
\end{enumerate}

\section{Numerical tests}
In this section, we consider the FSI-1 benchmark \cite{HrTu06b} 
(see also the books \cite{BuSc06,BuSc10} and our own former results \cite{RiWi10,Wi13_fsi_with_deal})
and the 2D-1 benchmark \cite{SchaeTu96}.
The drag value is taken as goal functional. As previously mentioned,
this paper is accompanied with a respective open-source 
implementation on github\footnote{\url{https://github.com/tommeswick/goal-oriented-fsi}} based on the finite element 
library deal.II \cite{dealII91,deal2020} and our previous fluid-structure interaction 
code \cite{Wi13_fsi_with_deal}, which is also available 
on github\footnote{\url{https://github.com/tommeswick/fsi}}.

\subsection{FSI-1 benchmark}
\begin{figure}[h]
\centering
{\includegraphics[width=7cm]{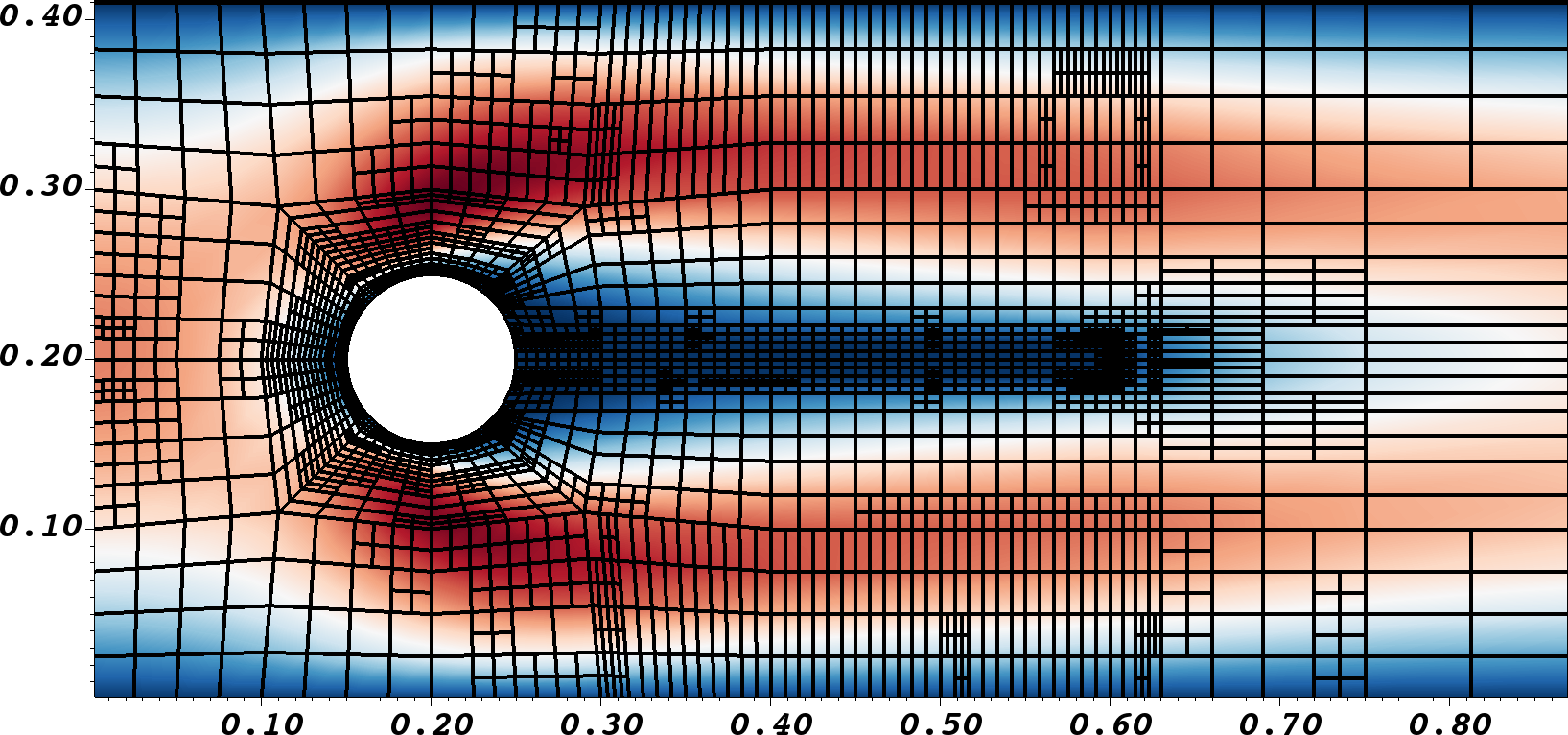}}
{\includegraphics[width=7cm]{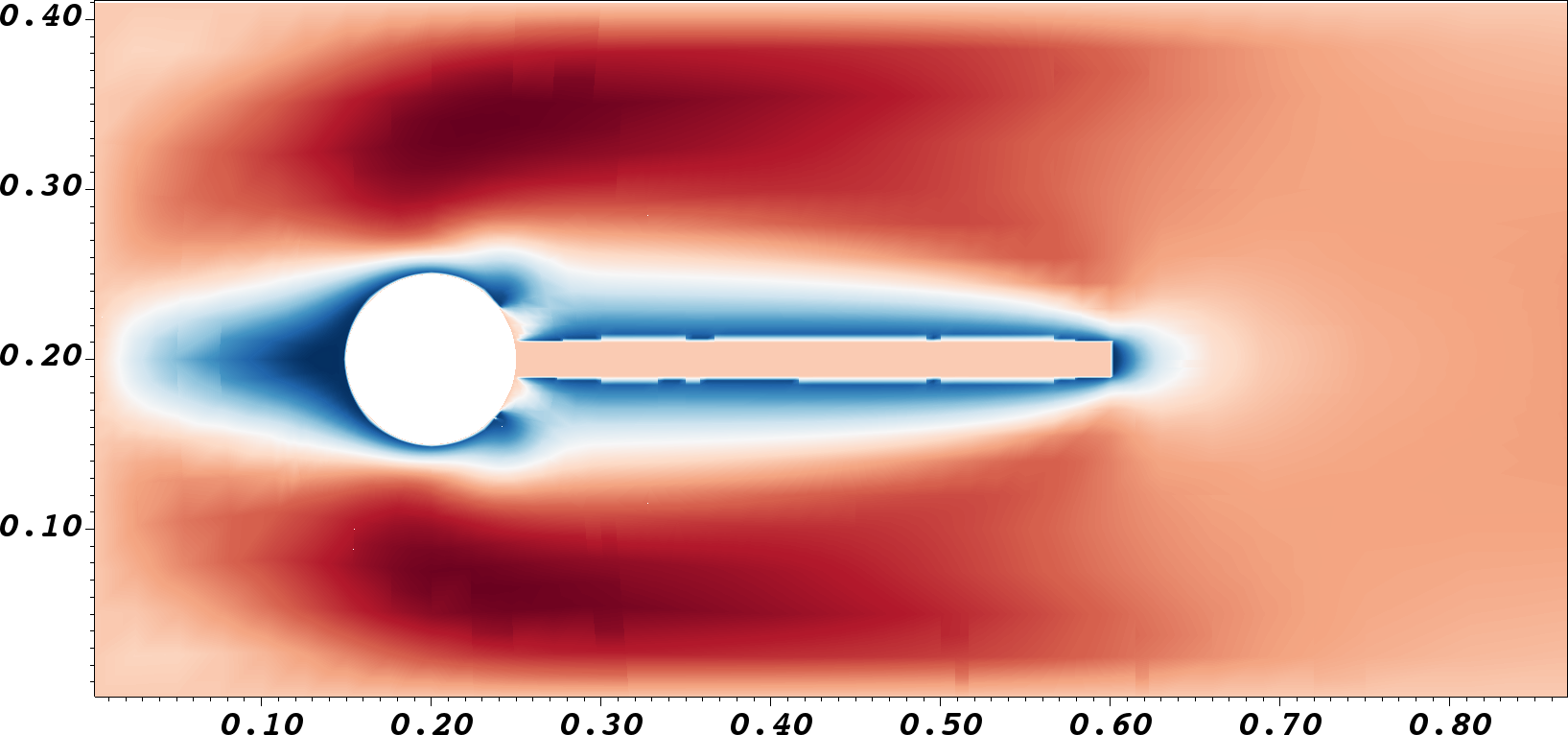}}
\caption{FSI-1 benchmark: Primal solution of $\hv_x$ and adjoint solution $\hz^{v_x}$. The adaptive mesh 
is displayed together with the primal solution (right).}
\label{pic_1}
\end{figure}

The configuration, all parameters, and reference values 
can be found in \cite{HrTu06b}. The reference value 
for computing the true error was computed on a five times refined mesh
and is $1.5370185576528707e+01$ (see also in the provided github code).
Our results from the file \verb|dwr_results.txt| are:
\begin{verbatim}
Dofs    True err     Est err      Est ind      Eff          Ind
13310   2.58e-01     1.43e-01     4.37e-01     5.54e-01     1.69e+00
20921   9.00e-02     4.75e-02     1.60e-01     5.28e-01     1.77e+00
37874   3.20e-02     1.09e-02     5.96e-02     3.40e-01     1.86e+00
68754   1.84e-02     4.57e-03     2.77e-02     2.48e-01     1.51e+00
\end{verbatim}
Furthermore, the terminal output yields
\begin{verbatim}
DisX  :     2.2656126465725842e-05
DisY  :     8.1965770448936843e-04
P-Diff:     1.4819455817646477e+02
P-front:    1.4819455817646477e+02
------------------
Face drag:      1.5351806985399641e+01
Face lift:      7.3933527637991259e-01
\end{verbatim}
where \verb|Face drag| represents the chosen goal functional. 
While the error reductions 
in the \verb|True err| $J(\hU) - J(\hU_h)$ 
and the estimated error $\eta$ are reasonable, the effectivity 
index \verb|Eff| has room for improvement. 
The indicator index \verb|Ind| (for the definition 
see \cite{RiWi15_dwr}) performs quite well. 
The main reason for the intermediate 
effectivity indices might be the accuracy of the reference value. Second, 
we notice that only the primal error part $\rho$ (Corollary \ref{coro_1}) 
was used. As shown 
in our recent studies for quasi-linear problems, the adjoint 
error part $\rho^*$ might play a crucial role in order to obtain nearly perfect effectivity 
indices for highly nonlinear problems \cite{EndtLaWi18}.
Graphical solutions of the primal solution, including the adaptively refined 
mesh, and the adjoint solution are displayed in Figure \ref{pic_1}.

\subsection{Adaptation to flow benchmark 2D-1}
The provided code can be adapted with minimal changes 
to the 1996 flow around cylinder benchmark 2D-1 \cite{SchaeTu96}.
In the \texttt{*.inp} file the material 
ids for solid must be set to $0$ (flow),
and the inflow and material parameters are adapted correspondingly. 
Of course, in this code,
the displacement variables are still computed despite that they are zero everywhere,
which increases the computational cost in comparison to a pure fluid flow code.
Due to the zero displacements $\hu=0$, the ALE mapping is the identity, 
yielding
$\hF=\hI$ and $det(\hF) = 1$. Consequently, there is no mesh deformation and 
the Navier-Stokes equations fully remain in Eulerian coordinates.
Here, extracting information from \verb|dwr_results.txt|, the findings 
for the drag value as goal functional are:
\begin{verbatim}
Dofs    True err    Est err      Est ind      Eff          Ind
1610    3.51e-01    2.97e-01     6.20e-01     8.44e-01     1.76e+00
2586    8.80e-02    7.27e-02     2.21e-01     8.26e-01     2.51e+00
4764    1.89e-02    1.54e-02     7.11e-02     8.14e-01     3.75e+00
10830   3.23e-03    2.95e-03     1.82e-02     9.13e-01     5.62e+00
\end{verbatim}
The pressure, drag (goal functional), and lift values are 
taken from the terminal output:
\begin{verbatim}
P-Diff:     1.1743527755157424e-01
P-front:    1.3213237901562136e-01
P-back:     1.4697101464047121e-02
------------------
Face drag:      5.5754969431700365e+00
Face lift:      1.0717678080199560e-02
\end{verbatim}
These values fit well with the reference values given in \cite{SchaeTu96}.
Moreover, we observe very stable effectivity indices, which indicate 
that the primal error estimator $\rho$ (Corollary \ref{coro_1}) is for incompressible
Navier-Stokes a sufficient choice. Indeed, using this 
part only, was already suggested in early work \cite{BeRa01,BrRi06}. 
Finally, we notice that extensions to multiple goal functionals for the 
2D-1 benchmark were undertaken in \cite{endtmayer2019hierarchical,EndtLaWi20}.

\section{Conclusions}
In this work, we developed and implemented  PU-DWR goal-oriented error control
and spatial mesh adaptivity for stationary fluid-structure interaction.
An important part is the open-source programming code published on github.
As numerical example, the FSI-1 benchmark is chosen. Therein,
mesh adaptivity performs as expected and also the error reductions 
in the true error and estimated error are good. However, the effectivity index
may be improved. Extensions of this work include inter alia 
the implementation of the adjoint error part $\rho^*$, local-higher order 
interpolations for the adjoint rather than using 
global-higher order finite elements,
parallel iterative/multigrid linear solvers within Newton's method, 
and a 3D implementation. The latter 
is implementation-wise not difficult with deal.II's 
dimension-independent programming, but 
the linear solver becomes really important.

\section{Acknowledgments}
This work is supported by the Deutsche Forschungsgemeinschaft (DFG) 
under Germany’s Excellence Strategy within 
the cluster of Excellence PhoenixD (EXC 2122, Project ID 390833453).


\end{document}